\documentclass{amsart}
\newtheorem{Theorem}{Theorem}[section]
\newtheorem{Definition}[Theorem]{Definition}
\newtheorem{Proposition}[Theorem]{Proposition}

\newtheorem{Lemma}[Theorem]{Lemma}
\newtheorem{Corollary}[Theorem]{Corollary}
\theoremstyle{remark}

\newtheorem{Example}[Theorem]{Example}

\def\il{\int\limits_}

\def\eps{\varepsilon}

\def\ovr{\overline}

\def\om{\omega}
\def\Om{\Omega}
\def\al{\alpha}

\def\dl{\delta}

\def\bd{\partial}

\def\si{\sigma}

\def\sm{\setminus}
\def\sbs{\subset}

\def\sea{\searrow}

\def\re{{\mathbf {Re\,}}}
\def\im{\mathbf {Im\,}}
\def\be{\begin{enumerate}}
\def\ee{\end{enumerate}}
\def\bT{\begin{Theorem}}
\def\eT{\end{Theorem}}
\def\bP{\begin{Proposition}}
\def\eP{\end{Proposition}}
\def\bD{\begin{Definition}}
\def\eD{\end{Definition}}
\def\bE{\begin{Example}}
\def\eE{\end{Example}}
\def\bL{\begin{Lemma}}
\def\eL{\end{Lemma}}
\def\bC{\begin{Corollary}}
\def\eC{\end{Corollary}}
\def\A{{\mathcal A}}

\def\aD{\mathbb D}

\def\aC{\mathbb C}

\def\aN{\mathbb N}

\def\E{{\mathcal E}}

\begin{document}
\title{Projective limits of Poletsky--Stessin Hardy spaces}
\author{Evgeny A. Poletsky}
\begin{abstract} In this paper we show that that on a strongly pseudoconvex domain $D$ the projective limit of all Poletsky--Stessin Hardy spaces $H^p_u(D)$, introduced in \cite{PS}, is isomorphic to the space $H^\infty(D)$ of bounded holomorphic functions on $D$ endowed with a special topology.
\par To prove this we show that Carath\'eodory balls lie in approach regions, establish a sharp inequality for the Monge--Amp\'ere mass of the envelope of plurisubharmonic exhaustion functions and use these facts to demonstrate that the intersection of all Poletsky--Stessin Hardy spaces $H^p_u(D)$ is $H^\infty(D)$.
\end{abstract}
\thanks{The author was partially supported by a grant from Simons Foundation.}
\keywords{Hardy spaces, pluripotential theory}
\subjclass[2000]{ Primary: 32A35; secondary: 32A70, 32U10}
\address{Department of Mathematics,  Syracuse University, \newline
215 Carnegie Hall, Syracuse, NY 13244} \email{eapolets@syr.edu}
\maketitle
\section{Introduction}
\par In \cite{PS} M. Stessin and the author introduced on a general hyperconvex domain $D$ the  spaces of holomorphic functions $H^p_u(D)$ as analogs of the classical Hardy spaces on the unit disk. This spaces are  parameterized by plurisubharmonic exhaustion functions $u$ of $D$. When $D$ is strictly pseudoconvex they all are the subsets of classical Hardy spaces $H^p(D)$ studied, for example, in \cite{S} and coincide with $H^p(D)$ when $u$ is a pluricomplex Green function.
\par Recently, M. Alan and N. Gogus in \cite{AG}, S. Sahin in \cite{Sa}, K. R. Shrestha in \cite{Sh1} and the latter with the author in \cite{PSh} showed that if $D$ is the unit disk $\aD$ these spaces form a subclass of weighted Hardy spaces studied, for example, in \cite{M} and \cite{BPST}. However, these subclass has special properties and, moreover, has no analogs in several variables. That is why we kept for it the name of Poletsky--Stessin Hardy spaces that is already used in these papers.
\par The parametrization of these spaces by plurisubharmonic exhaustion functions transforms this class into a projective system. In this paper we show that on a strongly pseudoconvex domain $D$ the projective limit of this system can be identified with the space $H^\infty(D)$ of bounded holomorphic functions on $D$ endowed with the projective topology. To prove this we construct for any unbounded holomorphic function $f$ a plurisubharmonic exhaustion function $u$ such that $f\not\in H^p_u(D)$. The construction is based on sharp estimates of the total Monge--Amp\'ere mass of the plurisubharmonic envelope of exhaustion functions (see Section \ref{S:mam}) and a placement of Carath\'eodory balls into Stein's approach regions in Section \ref{S:argb}.
\par We are grateful to M. Alan, N. Gogus, S. Sahin and K. R. Shrestha for stimulating discussions.
\section{Approach regions and balls}\label{S:argb}
\par Let $D$ be a bounded domain in $\aC^n$ with $C^2$ boundary. For $z_0\in\bd D$ we denote by $\nu_{z_0}$ the unit outward normal to $\bd D$ at $z_0$. Following E. Stein in \cite{S} for $\al>1$ we define the approach region $\A_D^\al(z_0)$ at $z_0$ as \[\A_D^\al(z_0)=\{z\in D:\, |(z-z_0)\cdot\nu_{z_0}|<\al\dl_D(z), |z-z_0|^2<\al\dl_D(z)\},\]
where $\dl_D(z)$ is the minimum of the distances from $z$ to $\bd D$ or to the tangent plane to $\bd D$ at $z_0$.
\par Recall that the Carath\'eodory function $c(z,w)$ on $D$ is defined as the supremum of $|f(z)|$ over all holomorphic functions $f$ on $D$ such that $f(w)=0$ and $|f|\le1$ on $D$. We define Carath\'eodory balls  centered at $w$ and of radius $r<1$ as the sets   $C_D(w,r)=\{z\in D:\,c(z,w)\le r\}$.
\par We will need the following result (see \cite[Theorem 2]{G}).
\bT\label{T:Gt}  Let $D$ be a strongly pseudoconvex domain in $\aC^n$ with $C^2$ boundary and $z_0\in\bd D$. Let $p$ be a peak function on $D$ at $z_0$, i.e., $p$ is continuous on $\ovr D$, holomorphic on $D$, $p(0)=1$ and $|p|<1$ elsewhere on $\ovr D$. Let $0<a<b<1$ and let $S(a)=\{z\in D:\,|p(z)|>a\}$. Choose any $\eta>0$. Then there exists  a positive constant $L=L(D,a,b,\eta)\ge1$ such that the following holds: given $f\in H^\infty(S(a))$, there exists $\hat f\in H^\infty(D)$ such that $\|\hat f\|_{H^\infty(D)}\le L\|f\|_{H^\infty(S(A))}$ and $\|f-\hat f\|_{H^\infty(S(b))}\le\eta\|f\|_{H^\infty(S(A))}$.
\eT
\bL\label{L:balls} Let $D$ be a strongly pseudoconvex domain in $\aC^n$ with $C^2$ boundary and $z_0\in\bd D$. For every $0< r<1$  there is $\al>0$ with the following property: for every neighborhood $U$ of $z_0$ there is $z\in D\cap U$ such that the Carath\'eodory  ball $C_D(z,r)$ lies in the approach region $\A_D^\al(z_0)$.
\eL
\begin{proof} We will prove this lemma in steps.
\par {\bf Step 1:} {\it The lemma holds when $D$ is the unit ball $B$ centered at the origin and $z_0=(1,0,\dots,0)$. One can take as $z$ any point $z=tz_0$, $0<t<1$, and $\al=20(1-r)^{-1}$.}
\par  Since $B$ has a transitive group of biholomorphisms, $C_B(v,r)=F(C_B(0,r))$, where $F$ is a biholomorphism of $B$ moving $0$ into $v$. Note that $C_B(0,r)$ is the ball of radius $r$ centered at the origin.
\par We let $v=(t,0,\dots,0)$, where $0<t<1$. If $z=(z_1,\dots,z_n)\in\aC^n$ then we set $z'=(z_2,\dots,z_n)$. The biholomorphism $(w_1,w')=F(z_1,z')$ moving 0 to $v$ is given by the formulas:
\[w_1=\frac{t+z_1}{1+tz_1}\text{ and }w'=(1-t^2)^{1/2}\frac{z'}{1+tz_1}.\]
Since for the ball the distance from a point in the ball to the boundary never exceeds the distance to the tangent plane $\dl_B(w_1,w')=1-(|w_1|^2+|w'|^2)^{1/2}$.
If $(z_1,z')\in C_B(0,r)$ then
\[\dl_B(w_1,w')\ge \frac12(1-|w_1|^2-|w'|^2)\ge\frac{(1-t^2)(1-r^2)}{2|1+tz_1|^2}\ge\frac{(1-t)(1-r)}{2|1+tz_1|^2}. \]
Since $\nu_{z_0}=(1,0')$ for $(w_1,w_2)\in C_B(v,r)$ we have
\[|(w-z_0)\cdot\nu_{z_0}|=|1-w_1|=\frac{(1-t)|1-z_1|}{|1+tz_1|}\le\frac{4(1-t)}{|1+tz_1|^2}\]
and
\[|w-z_0|^2=|w'|^2+|1-w_1|^2=\frac{(1-t^2)|z'|^2+(1-t)^2|1-z_1|^2}{|1+tz_1|^2}\le
\frac{10(1-t)}{|1+tz_1|^2}.\]
\par Therefore, for every $0<t<1$ the Carath\'eodory ball $C_B((t,0'),r)$ lies in the approach region $\A_B^\al(z_0)$, when $\al=20(1-r)^{-1}$ and this ends Step 1.
\par {\bf Step 2:} {\it Let $0\le t<1$, $z_0=(1,0,\dots,0)$, $B_{-t}=\{z\in\aC^n:\,|z+tz_0|<1+t\}$ and $B_t=\{z\in\aC^n:\,|z-tz_0|<1-t\}$. Then  $\A^\al_{B_{-t}}(z_0)\sbs\A^{4\al}_{B_t}(z_0)$ when $0<t<(8\al)^{-1}$.}
\par If $z=(z_1,\dots,z_n)\in B_{-t}$ and $x=\re z_1$, then
$\dl_{B_{-t}}(z)=1+t-|z+tz_0|$ and $\dl_{B_t}(z)=1-t-|z-tz_0|$. Direct calculations show that
\[(1+t+|z+tz_0|)\dl_{B_{-t}}(z)=(1-t+|z-tz_0|)\dl_{B_t}(z)+4t(1-x).\] Thus $\dl_{B_{-t}}(z)\le2\dl_{B_t}(z)+4t(1-x)$. But $1-x<\al\dl_{B_{-t}}(z)$. Hence $\dl_{B_{-t}}(z)\le 2(1-4t\al)^{-1}\dl_{B_t}(z)$. If $0<t<(8\al)^{-1}$ then $\dl_{B_{-t}}(z)\le 4\dl_{B_t}(z)$. So if $z\in\A_{B_{-t}}^\al(z_0)$ then $z\in\A_{B_{-t}}^{4\al}(z_0)$.
\par {\bf Step 3:} {\it Let $p$ be a peak function at $z_0$. If the lemma holds for some $S(a)=\{z\in D:\,|p(z)|>a\}$ then it holds for $D$.}
\par The function $\dl(z)$ in the definition of approach regions is the same whether we take it with respect to $D$ or $S(a)$ when $z$ is sufficiently close to $z_0$. So we can take $b_0$, $a<b_0<1$ so that the intersections of approach regions with respect to $D$ or $S(a)$ coincide in $S(b_0)$.
\par Fix some positive  $r<1$ and let $r'=r+(1-r)/2$. We take $\eps,\eta>0$ such that
\[(1+2\eta)^{-1}(1-\eps)(r'-2\eta)>r.\]
Let $L=L(D,a,b_0,\eta)$. We take an integer $m$ such that $b_0^mL<1$  and a number $b$ between $b_0$ and 1 such that $b^m>1-\eps$. There is $c$, $b<c<1$, such that the Carath\'eodory balls $C_D(w,r')\sbs S(b)$ when $w\in S(c)$. Indeed, if $z_0\in C(w,r')$ and $f$ is a conformal mapping $f$ of the unit disk onto itself such that $f(p(w))=0$, then $|f(p(z_0)|\le r'$. Direct calculations show that if $|p(w)|>(b+r')/(1+br')$ then $|p(z_0)|>b$.
\par Since the lemma holds on $S(a)$ we can find $\al$ and $w_0\in S(c)$ such that for every point $w\not\in\A_D^\al(z_0)\cap S(b)$ there is a holomorphic function $f$ on $S(a)$ such that $|f|<1$ on $S(a)$, $f(w_0)=0$ and $|f(w)|>r'$. By Theorem \ref{T:Gt} there is a function $\hat f\in H^\infty(D)$ such that $\|\hat f\|_{H^\infty(D)}\le L$ and $\|f-\hat f\|_{H^\infty(S(b_0))}\le\eta$.
\par Let $g=(1+2\eta)^{-1}p^m(\hat f-\hat f(w_0)$. If $z\in D\sm S(b_0)$ then $|g(z)|\le b_0^mL\le 1$.
If $z\in S(b_0)$ then $|g(z)|<(1+2\eta)^{-1}(1+2\eta)=1$. Hence $|g|<1$ on $D$. Now
\[|g(w)|\ge(1+2\eta)^{-1}b^m(r'-2\eta)>(1+2\eta)^{-1}(1-\eps)(r'-2\eta)>r.\]
Hence $w\not\in C_D(w_0,r)$ and $C_D(w_0,r)\sbs\A_D^\al(z_0)\cap S(b)$. This ends Step 3.
\par We take a plurisubharmonic function $\phi\in C^2(\ovr D)$ defining $D$ such that $\nabla\phi\ne0$ on $\bd D$.
Let \[L_{z_0}(z)=\sum_{i,j=1}^n\phi_{z_i,z_j}(z_0)(z_i-(z_0)_i)(z_j-(z_0)_j)\] and  \[H_{z_0}(z)=\sum_{i,j=1}^n\phi_{z_i,\ovr z_j}(z_0)(z_i-(z_0)_i)(\ovr z_j-(\ovr z_0)_j).\]
The Taylor expansion of $\phi$ at $z_0$ is
\[\phi(z)=2\re (\nabla\phi(z_0),z-z_0) +\re L_{z_0}(z)+\frac12H_{z_0}(z)+o(\|z-z_0\|^2).\]
\par {\bf Step 4:} {\it The lemma holds when $z_0=(1,0,\dots,0)$ and the Taylor expansion of $\phi$ at $z_0$ is
\[\phi(z)=-2(1-x)+|z-z_0|^2 +o(|z-z_0|^2).\]}
\par We take $\al=20(1-r)^{-1}$ and $t=(16\al)^{-1}$. By Step 1
$C_B(sz_0,r)\sbs\A_B^\al(z_0)$ for any $0<s<1$. The dilation $d(z)=(1+t)z-tz_0$ moves $B$ onto $B_{-t}$ and $C_B(sz_0,r)$ onto $C_{B_{-t}}(s'z_0,r)$, $s'=(1+t)s-t$. If $z\in\A_B^\al(z_0)$ then
\[|d(z)-z_0|^2=(1+t)^2|z-z_0|^2<(1+t)^2\al\dl_B(z)=(1+t)\al\dl_{B_{-t}}(z)<
2\al\dl_{B_{-t}}(z),\]
while $|(d(z)-z_0)\cdot\nu_{z_0}|=\al\dl_{B_{-t}}(z)$.
Thus $d$ moves $\A_B^\al(z_0)$ into $\A_{B_{-t}}^{2\al}(z_0)$ and we see that $C_{B_{-t}}(sz_0,r)\sbs\A_{B_{-t}}^{2\al}(z_0)$ for any $0<s<1$.
\par There is $x_0<1$ such that if $\Om=\{z\in D:\,\re z_1>x_0\}$, $B'_t=B_t\cap\Om$, $B'_{-t}=B_{-t}\cap\Om$ and $D'$ is the connected component of $D\cap\Om$ containing $z_0$, then $B'_t\sbs D'\sbs B'_{-t}$. Hence by Step 2
\[C_{D'}(sz_0,r)\sbs C_{B'_{-t}}(sz_0,r)\sbs\A_{B_{-t}}^{2\al}(z_0)\sbs\A_{B_t}^{8\al}(z_0)\sbs\A_{D'}^{8\al}(z_0)\] when $s$ is sufficiently close to $1$. By Step 3 the statement holds.
\par {\bf Step 5:} {\it General case.} There is (see Lemma 5 and Proposition 2 in \cite{G}) a quadratic transformation $F$ of $\aC^n$, biholomorphic in a neighborhood $U$ of $z_0$, that moves $D$ into a domain where the Taylor expansion of $\phi$ at $\phi(z_0)$ has the form
\[\phi(z)=-2\re z_1+\sum_{j=1}^n|z-z_0|^2+o(|z-z_0|^2).\]
Since the image and the preimage of approach regions under the mapping $F$ will lie in corresponding approach regions near the boundary by Steps 3 and 4 we get our lemma.
\end{proof}
\section{The Monge--Amp\'ere mass of envelopes}\label{S:mam}
\par A domain $D\sbs\aC^n$ is {\it hyperconvex} if there is a continuous function $u$ on $\ovr D$ equal to zero on $\bd D$ and negative and plurisubharmonic on $D$ and it is {\it strongly hyperconvex} if $u$ extends as a continuous plurisubharmonic function to a neighborhood of $\ovr D$.  We denote by $\E(D)$ the set of all continuous functions $u$ on $\ovr D$ equal to zero on $\bd D$ and negative and plurisubharmonic on $D$. We assume that such functions can take $-\infty$ as their value.
\par A {\it pluriregular condensor }
$K=(K_1,\dots,K_m,\sigma_1,\dots,\sigma_m)$ is a system of pluriregular compact sets \[K_m\subset K_{m-1}\subset\dots\subset K_1\subset D\subset\overline D=K_0\] and numbers
$\sigma_m<\sigma_{m-1}<\dots<\sigma_1<\sigma_0=0$ such that there is a continuous plurisubharmonic function $\omega(z)=\omega(z,K,D)$ on $D$ with zero boundary values,
$K_i=\{\omega\le\sigma_i\}$ and $\omega$ is maximal on $D_{\si_{i-1}}\sm K_i$ for all $1\le i\le m$, where $D_\si=\{z\in D:\,\om(z)<\si\}$ (see \cite{P} for more details). We will call this function {\it the relative extremal function} of the condensor $K$ in $D$. Of course, not every choice of sets $K_i$ and numbers $\sigma_i$ can be realized as a condensor. But if $u$ is a continuous negative plurisubharmonic function on $D$ and sets $K_i=\{u\le\sigma_i\}$ are pluriregular,
then $K$ has a continuous relative extremal function.
\par The following lemma was proved in \cite[Lemma 4.2]{P}.
\begin{Lemma}\label{L:ca} Let $K=(K_1,\dots,K_m,\sigma_1,\dots,\sigma_m)$ be a pluriregular condensor in $D$. There is a sequence of pluricomplex multipole Green functions $g_j$ converging to $\omega(z)=\omega(z,K,D)$ uniformly on
compacta in $D_{\sigma_{i-1}}\setminus K_i$, $1\le i\le m$.
Moreover, if $\psi$ is a continuous function on ${\mathbb  R}$, then
\[\lim_{j\to\infty}\int\limits_ D\psi(\omega(z))(dd^cg_j)^n=
\int\limits_ D\psi(\omega(z))(dd^c\omega)^n.\]\end{Lemma}
\par The following lemma is a slight but important elaboration of the previous result.
\bL\label{L:cae} Let $K=(K_1,\dots,K_m,\sigma_1,\dots,\sigma_m)$ be a pluriregular condensor in $D$. There is a sequence of pluricomplex multipole Green functions $g_j(z)<\om(z)=\omega(z,K,D)$ on $D$ converging to $u(z)$ uniformly on compacta in $D_{\sigma_{i-1}}\setminus K_i$, $1\le i\le m$. Moreover, if $\psi$ is a continuous function on ${\mathbb  R}$, then
\[\lim_{j\to\infty}\int\limits_ D\psi(\omega(z))(dd^cg_j)^n=
\int\limits_ D\psi(\omega(z))(dd^c\omega)^n.\]
\eL
\begin{proof} By Lemma \ref{L:ca} there is a sequence of Green functions $h_j$ on $D$ converging to $u$ uniformly on
compacta in $D_{\sigma_{i-1}}\setminus K_i$, $1\le i\le m$. Moreover, if $\psi$ is a continuous function on ${\mathbb  R}$, then
\[\lim_{j\to\infty}\int\limits_ D\psi(\om(z))(dd^ch_j)^n=
\int\limits_ D\psi(\om(z))(dd^c\om)^n.\]
\par Let us choose a decreasing sequence of numbers $\al_k>1$ converging to 1. Define $\si'_{ik}=\al_k^{-1}\si_i$ and $\si''_{ik}=\al_k\si_i$. There is $k_0$ such that for all $k>k_0$ and $i=1,\dots,m$ we have
\[\si_{i+1}<\si_{ik}''<\si_i<\si'_{ik}<\si_{i-1}.\]
\par For any such $k$ there is $j_k>k$ such that $\al_kh_{j_k}<\si'_{ik}$ on $\bd D_{\si'_{ik}}$ and $\al_kh_{j_k}<\si''_{ik}$ on $\bd D_{\si''_{ik}}$. By the maximum principle $\al_kh_{j_k}<\si'_{ik}$ on $D_{\si'_{ik}}\sm D_{\si''_{ik}}$. Hence
\[\al^3_kh_{j_k}<\al_k^2\si'_{ik}=\al_k\si_i=\si''_{ik}\] on $D_{\si'_{ik}}\sm D_{\si''_{ik}}$. So if $g_k=\al_k^3h_{j_k}$ then $g_k<\om$ on $D_{\si'_{ik}}\sm D_{\si''_{ik}}$ for all $i=1,\dots,m$. By the maximality of $\om$ on $D_{\si_i}\sm\ovr D_{\si_{i+1}}$ and we see that $g_k<\om$ on $D$. Clearly,
\[\lim_{j\to\infty}\int\limits_ D\psi(\omega(z))(dd^cg_j)^n=
\int\limits_ D\psi(\omega(z))(dd^c\omega)^n.\]
\end{proof}
\par Given a continuous function $\phi$ on $D$ we denote by $E\phi$ the plurisubharmonic envelope of $\phi$, i.e., the maximal plurisubharmonic function on $D$ less or equal to $\phi$. Such a function exists due to the continuity of $\phi$. By \cite[Lemma 1]{W} if $\phi<0$ on $D$ and $\lim_{z\to\bd D}\hat\phi(z)=0$, the $E\phi$ is continuous on $D$.  For an at most countable sequence of functions $\{u_j\}\sbs\E$ we denote by $E\{u_j\}$ the envelope of $\min\{u_j\}$.
\bT\label{T:mae} If $D$ is a strongly hyperconvex domain and continuous plurisubharmonic functions $\{u_j\}\sbs\E(D)$, then $\sum MA(u_j)\ge MA(E\{u_j\})$, where
\[MA(u)=\int\limits_ D(dd^cu)^n.\]
\eT
\begin{proof} First, we prove this theorem for two functions $u$ and $v$. Since $E(u,v)\ge u+v$ we see that $E(u,v)\in\E(D)$. We may assume that functions $u,v\in\E(D)$ are bounded and of the finite Monge--Amp\'ere mass. If the former does not hold then we replace $u$ and $v$ with $u_k=\max\{u,-k\}$ and $v_k=\max\{v,-k\}$ respectively and use the fact that $MA(u_k)\to MA(u)$ for a decreasing sequence $\{u_k\}$ and $E(u_k,v_k)\sea E(u,v)$. If the latter is not true then the statement is evident.
\par If $K$ and $L$ are pluriregular condensors in $D$, $u(z)=\om(z,K,D)$ and $v(z)=\om(z,L,D)$, then by Lemma \ref{L:cae} there are  sequences of pluricomplex multipole Green functions $\{g_j(z)<u(z)\}$ and $\{h_j(z)<v(z)\}$ on $D$ such that
\[\lim_{j\to\infty}MA(g_j)=MA(u)\] and
\[\lim_{j\to\infty}MA(h_j)=MA(v).\]
\par Clearly, $E(u,v)>E(g_j,h_j)$ and by the Comparison Principle
$MA(E(u,v))\le MA(E(g_j,h_j))$.
But $E(g_j,h_j)$ is a pluricomplex multipole Green function with poles at poles of $g_j$ and $h_j$ and weights equal to the maximum of weights $g_j$ or $h_j$ at a pole. Hence
$MA(E(g_j,h_j))\le MA(g_j)+MA(h_j)$ and our theorem holds in this case.
\par In the next step we prove the theorem for functions $u,v\in\E(D)$ for which there is an open set $D'\subset\subset D$ such that $\partial D'$ is a smooth hypersurface, $u$ and $v$ are equal to $\sigma_1<0$ on $\partial D'$, maximal on $D\setminus\overline D'$ and are of class $C^2$ on $D'$.
\par For this we will construct an inductive sequence of pluriregular condensors $K_j$ and $L_j$ such that the sequences of functions $u_j(z)=\om(z,K_j,D)$ and $v_j(z)=\om(z,L_j,D)$ are decreasing and converging to $u$ and $v$ respectively. Then $E(u,v)$ is the limit of the decreasing sequence of $E(u_j,v_j)$ and, consequently, the theorem holds in this case.
\par We let $K_0=(\ovr D',\si_1)$. If $K_j=(K_{1j}=\ovr D',K_{2j}\dots,K_{i_jj},\si_{1j}=\si_1,\si_{2j},\dots,\si_{m_jj})$ has been constructed, then by Sard's theorem for every $1\le i\le m_j-1$ we can find numbers \[\si_{j,i+1}=\dl_{l_{ij}}<\dl_{l_{ij}-1}<\dots<\dl_1<\dl_0=\sigma_{ij}\]
such that $\dl_l-\dl_{l+1}<1/j$ and the function $u$ is not degenerate on $\{u=\dl_l\}$, $1\le l\le l_{ij}$. For $i=m_j$ we select numbers $\dl_l$ as before between $\si_{i_jj}$ and the minimum of $u$ on $D$.
\par Since the hypersurfaces of $\{u=\dl_{l_{ij}}\}$ are smooth, the compact sets  $K_{l_{ij}}=\{u\le\dl_{l_{ij}}\}$ are pluriregular. We relabel the numbers $\si_{kj}$ and $\dl_{l_{ij}}$ and compact sets $K_{l_{ij}}$ as $\si_{i,j+1}$ and $K_{i,j+1}$ respectively arranging them in the right order and define a pluriregular condensor
\[K_{j+1}=(\ovr D',K_{2,j+1}\dots,K_{m_{j+1},j+1},\si_{1,j+1},\dots,\si_{m_{j+1},j+1}).\] We denote $\om(z,K_{j+1},D)$ by $u_{j+1}(z)$.
\par Since the functions $u_j$ are maximal on $K_{ij}^o\sm K_{i+1,j}$ and are equal to $u$ on $\bd K_{ij}$ we see that $u_j\ge u$ on $D$ and $u_j\ge u_{j+1}$ on $D$. Hence the sequence of $u_j$ is decreasing and, clearly, converging to $u$. Since a similar construction works for the function $v$ too, our theorem holds in this case.
\par For the general case we suppose $D=\{z\in\aC^n:\,\phi(z)<0\}$, where $\phi$ is a continuous plurisubharmonic function defined on a neighborhood $V$ of $\overline D$ and its restriction to $D$ is in $\E(D)$.  The sequence of plurisubharmonic functions $u_k$ on $V$ equal to $\max\{u,k\phi\}$ on $D$ and to $k\phi$ on $V\setminus D$ is decreasing on $D$ and converges to $u$ uniformly on $\overline D$. In particular, the total Monge--Amp\'ere masses of $u_k$ on $D$ converge to this mass of $u$. Hence if we prove our theorem for continuous functions that admit a continuous plurisubharmonic extension to $V$, then  we prove it for all functions in $\E(D)$.
\par If $u$ is such a function then there is a decreasing sequence of plurisubharmonic functions $u_k$ on some domain $G$ containing $\ovr D$ that belong to $C^\infty(G)$ (see \cite[Theorem 2.9.2]{K}) and converge to $u$ uniformly on $\overline D$. Let $\eps_k=\sup_{z\in\bd D}u_k(z)$.
\par Let us choose a sequence of numbers $\sigma_{1k}<0$ converging to 0 such that for all $k$ the set $\{u_k=\sigma_{1k}\}$ is a smooth hypersurface compactly belonging to $D$. We define $u'_k$ as a function that is equal to $u_k-\eps_k$ on the set $W_k=\{u_k\le\sigma_{1k}\}$, to 0 on $\partial D$ and to be maximal on $D\setminus W_k$. These functions uniformly converge to $u$ and they are plurisubharmonic because $u'_k\ge u_k-\eps_k$ on $D\sm W_k$. Hence the total Monge--Amp\'ere masses of $u'_k$ on $D$ converge to this mass of $u$. Since
for functions like this our theorem is already proved, it is proved for all $u\in\E(D)$.
\par For finitely many functions $u_1,\dots,u_k$ the result follows immediately by induction: the envelope $v_k$ of $\min\{u_1,\dots,u_k\}$ is equal to the envelope of \[\min\{\min\{u_1,\dots,u_{k-1}\},u_k\}.\] For the infinite case we note that $E(\{u_j\})$ is the limit of the decreasing sequence of $v_k$ and the inequality follows from the classical result of Bedford and Taylor.
\end{proof}
\par This result is sharp. If $D$ is hyperconvex and $W=\{w_1,\dots,w_k\}\sbs D$, then the {\it pluricomplex Green
function with poles at the set $W$}  is a unique unction $g(z,W)\in\E(D)$ such that
$(dd^cg(z,W))^n=\sum_{j=1}^k(2\pi)^n\dl_{w_j}$, $|g_D(z,W)-\sum_{j=1}^k\log|z-w_j||$ is
bounded on $D$ and $g(z,W)$ is maximal outside $W$, i.e., $(dd^cg)^n=0$ on $D\sm W$.
\par If $u$ and $v$ are two pluricomplex Green functions with non-overlapping poles, then $E(u,v)$ is the pluricomplex Green functions whose set of poles is the union of poles of $u$ and $v$. Hence we have an equality in Theorem \ref{T:mae}.
\par We finish this section with the following observation. Let $\E_1(D)$ be the set of all $u\in\E(D)$ such that $MA(u)=1$.
\bC\label{C:mae} If $u,v\in\E_1(D)$ then $MA(E(u,v))\le 2$.\eC

\section{Poletsky--Stessin  Hardy spaces}
\par Let $D$ be a hyperconvex domain in ${\mathbb  C}^n$ and $u\in\E(D)$. Following \cite{D} we set $B_u(r)=\{z\in D:\,u(z)< r\}$ and $S_u(r)=\{z\in D:\,u(z)= r\}$. Let $$\mu_{u,r}=(dd^cu_r)^n-\chi_{D\sm B_u(r)}(dd^cu)^n,$$
where $u_r=\max\{u,r\}$. The measure $\mu_{u,r}$ is nonnegative and
supported by $S_u(r)$. In \cite[Theorem 1.7]{D} Demailly had proved the
following fundamental Lelong--Jensen formula.
\bT\label{T:ljf} For all $r<0$ and every plurisubharmonic function
$\phi$ on $D$
$$\mu_{u,r}(\phi)=\il D\phi\mu_{u,r}$$ is finite and
\begin{equation}\label{e:ljf}
\mu_{u,r}(\phi)-\il {B_u(r)}\phi(dd^cu)^n=\il
{B_u(r)}(r-u)dd^c\phi\wedge(dd^cu)^{n-1}.
\end{equation}
\eT
\par The last integral in this formula can be equal to $\infty$.
Then the integral in the left side is equal to $-\infty$. This
cannot happen if $\phi\ge0$.
\par The function
\[\Phi(r)=\il{B_u(r)}(r-u)dd^c\phi\wedge(dd^cu)^{n-1}\] is, evidently, increasing  and it follows that the function $\mu_{u,r}(\phi)$ is increasing and continuous from the left.
\par  As in \cite{PS} for $p\ge1$ we define the Hardy space $H^p_u(D)$ as the set of all holomorphic functions $f$ on $D$ such that
\[\limsup_{r\to0^-}\mu_{u,r}(|f|^p)<\infty.\]
\par Since $\mu_{u,r}(|f|^p)$ is an increasing function of $r$ for all $r<0$, we can replace $\limsup$ in the definition of this space by $\lim$. So we can introduce the norm on $H^p_u(D)$ as
\[\|f\|^p_{u,p}=\lim_{r\to0^-}\mu_{u,r}(|f|^p)=\il {D}|f|^p(dd^cu)^n-\il
{D}udd^c|f|^p\wedge(dd^cu)^{n-1}.\]
It was shown (see \cite[Theorem 4.1]{PS}) that the spaces $H^p_u(D)$ are Banach for $p\ge1$. \par The following theorem which is a direct consequence of \cite[Corollary 3.2]{PS} shows that faster decaying near the boundary of $D$ exhausting functions determine dominating norms.
\bT\label{T:it} Let $u$ and  $v$ be continuous plurisubharmonic
exhaustion functions on $D$ and let $F$ be a compact set in $D$
such that $bv(z)\le u(z)$ for some constant $b>0$ and all $z\in
D\sm F$. Then $H^p_v(D)\sbs H^p_u(D)$ and $\|f\|_{u,p}\le
b^{n/p}\|f\|_{v,p}$.
\eT
\par  Let $u=(u_1,\dots,u_k)\in\E^k_1$. Let $H^p_u(D)$ be the direct product $H^p_{u_1}(D)\times\cdots\times H^p_{u_k}(D)$ with the norm
\[\|(f_1,\dots,f_k)\|_{u,p}=\sum_{j=1}^k\|f_j\|_{u_j,p}.\]  We denote by $B_{u,p}(r)$ the open ball of radius $r$ centered at the origin of $H^p_u$.
\par If $p=2$ then we introduce a sesqui-linear form on $H^2_u(D)$ as
\[(f,g)_u=\lim_{r\to 0^-}\sum_{j=1}^k\il{S_u(r)}f_j\ovr g_j\,d\mu_{u,r}.\]
Since $2\re f\ovr g=|f+g|^2-|f|^2-|g|^2$ and $2\im f\ovr g=|f+ig|^2-|f|^2-|g|^2$, the existence of the limit follows. By the H\"older inequality $|(f,g)_u|^2\le(f,f)_u(g,g)_u<\infty$.
It follows that a continuous non-negative sesqui-linear form $(f,g)$ is well defined on $H^2_u(D)$ and makes this space a Hilbert space.
\par  The norm of $f=(f_1,\dots,f_k)\in(H^\infty(D))^k$ will be defined as
\[\|f\|_{\infty}=\sum_{j=1}^k\|f_j\|_\infty\] and let $B^k_\infty(r)$ be the open ball of radius $r$ centered at the origin of $(H^\infty)^k$. If $f,g\in H^p_u(D)$ and $|f|\le|g|$ on $D$, then $\mu_{u,r}(|f|^p)\le \mu_{u,r}(|g|^p)$. Hence $\|f\|_{u,p}\le\|g\|_{u,p}$ and we see that $B^k_\infty(r)\sbs B_{u,p}(r)$ when $u\in\E_1^k$.
\par If $u=(u_1,\dots,u_k)$ and $v=(v_1,\dots,v_k)$ are in $\E^k_1(D)$ then we say that  $u\succeq v$ if there is a constant $c>0$ and a compact set $F\sbs D$ such that $cu_j\le v_j$ on $D\sm F$. In this case $H^p_u(D)\sbs H^p_{v,k}(D)$ and there is a constant $a>0$ such that $\|f\|_{v,p}\le a\|f\|_{u,p}$.
\bP\label{P:cs} Let $u,v\in\E^k_1(D)$ and $u\succeq v$. Then:\be
\item If $A\sbs H^p_v(D)$ is closed in $H^p_v(D)$ then $A\cap H^p_u$ is closed in $H^p_u(D)$;
\item the closed balls $\ovr B_{u,p}(R)$ in $H^p_u(D)$ of radius $R$ are closed in $H^p_v(D)$;
\item if $A\sbs H^2_u(D)$ is a closed convex bounded set, then $A$ is a closed bounded set in $H^2_v(D)$.
    \ee
\eP
\begin{proof} (1) Indeed, if a sequence $\{f_j\}\sbs A\cap H^p_u(D)$ and $f_j\to f$ in $H^p_u(D)$, then  $\|f_j-f\|_{v,p}\le c\|f_j-f\|_{u,p}$. Hence $f_j\to f$ in $H^p_v(D)$ and $f\in A$.
\par (2) Let $\{f_j=(f_{j1},\dots,f_{jk})\}$ be a sequence in $\ovr B_{u,p}(R)$ converging in $H^p_v(D)$ to  $g=(g_1,\dots,g_k)$. Then the functions $f_{jm}$ converge to $g_m$ in $H^p_{v_m}(D)$. Since  the integrals $\mu_{u_m,r}(|f_{jm}|^p)$ are increasing in $r$ we see that $\mu_{u_m,r}(|f_{jm}|^p)\le \|f_{jm}\|^p_{u_m,p}$. By Theorem 3.6 from \cite{PS} $\{f_{jm}\}$ is a Cauchy sequence in the uniform metric on any compact set in $D$. Hence, for any $r<0$
\[\mu_{u_m,r}(|g_m|^p)=\lim_{j\to\infty}\mu_{u_m,r}(|f_{jm}|^p)\le
\lim_{j\to\infty}\|f_{jm}\|^p_{u_m,p}.\]
Consequently,
\[\|g\|_{u,p}\le \lim_{j\to\infty}\|f_j\|_{u,p}\le R\]
and we see that $f\in\ovr  B_{u,p}(R)$.
\par (3) The fact that $A$ is bounded in $H^2_v(D)$ follows from Theorem \ref{T:it}. Let $\{f_j=(f_{j1},\dots,f_{jk})\}$ be a sequence in $A$ converging in $H^2_v(D)$ to  $g=(g_1,\dots,g_k)$. Then the functions $f_{jm}$ converge to $g_m$ in $H^2_v(D)$. As it was observed in part (2) $\{f_{jm}\}$ is a Cauchy sequence in the uniform metric on any compact set in $D$.
\par Since Hilbert spaces are reflexive, the closed balls are weakly compact. Since $A$ is convex and closed it is weakly closed in $H^2_{u,k}(D)$. Hence there is a subsequence $\{f_{j_k}\}$ weakly converging to $h=(h_1,\dots,h_k)\in A$.
\par If $D$ is a hyperconvex domain, $w_0\in D$ and $\om(z,w_0)$ is the pluricomplex Green function with pole in $w_0$, then for any $u\in\E$ there is a constant $c>0$ such that $cu\le\om$ near $\bd D$. Hence $H^p_u(D)\sbs H^p_\om(D)$ and $\|f\|_{\om,p}\le c^{n/p}\|f\|_{u,p}$. By formula (3.2) in \cite{PS}
\[(2\pi)^n|f(w)|^p\le\|f\|_{p,\om}\le c^{n/p}\|f\|_{p,u}\] when $f\in H^p_u(D)$. Hence point evaluations are continuous functionals on $H^p_u(D)$. Thus $h_m=g_m$.
\end{proof}
\par The following result was proved by K. R. Shrestha in \cite{Sh2} when $D$ is the unit disk.
\bT\label{T:ib} If $\ovr B=\cap_{u\in\E_1^k}\ovr B_{u,p}(R)$ then $\ovr B=\ovr B_\infty(R)$.
\eT
\begin{proof} Suppose that $f=(f_1,\dots,f_k)\in\ovr B$ and $\sum\|f_j\|_\infty>R$. We fix an $\eps>0$ and find points $w_1,\dots,w_k\in D$ such that $|f_j(w_j)|\ge\|f_j\|_\infty-\eps$. Let $u_j(z)=g(z,w_j)$. Then $\|f_j\|_{u_j,p}\ge|f_j(w_j)|$ and $\|f\|_{u,p}\ge\|f\|_\infty-k\eps$, where $u=(u_1,\dots,u_k)$. Since $\eps$ is arbitrary we come to a contradiction.
\par If $f\in\ovr B_\infty(R)$ then $\|f_j\|_{u,p}\le\|f_j\|_\infty$ for any $u\in\E_1^k$. Hence $f\in \ovr B$.
\end{proof}
\par The following result gives some chances for a reduction of $H^\infty$ problems to $H^2$ problems.  When $D$ is the unit disk it was proved in \cite{PSh} by K. R. Shrestha and the author for any $p>1$ and without any conditions.
\bT \label{T:int1} Let $D$ be a strongly pseudoconvex domain, $w_0\in D$, $R>0$ and $g=g(z,w_0)$.  Let $A \sbs (H^2_g(D))^k$ be a closed convex set. Then $A\cap \ovr B_\infty(R) \neq \emptyset$ if and only if for any number of functions $u_1,\dots,u_m\in\E^k_1(D)$ the set $A\cap \ovr B_{u_1,2}(R)\cap\cdots\cap \ovr B_{u_m,2}(R) \neq \emptyset$. \eT
\begin{proof} By Proposition \ref{P:cs}(1) for $u\in\E^k_1(D)$ the set $A_u=A\cap \ovr B_{u,2}(R)$ is closed in $H^2_u(D)$. Since it is convex and bounded by Proposition \ref{P:cs}(3) it is closed and bounded in $(H^2_g(D))^k$. Since it is convex it is weakly closed in $(H^2_g(D))^k$ and, consequently, weakly compact.
\par Since any finite number of the sets $A_{u_1},\dots A_{u_m}$ have the non-empty intersection we see that $\cap_{u\in\E_1^k}A_u \neq \emptyset$. By Theorem \ref{T:ib} the latter set is equal to $A\cap B_\infty(R)$.
\par If $f\in A\cap \ovr B_\infty(R)$ then $\|f\|_{u,2}\le R$ and the theorem follows.
\end{proof}
\section{Projective limits of Poletsky--Stessin Hardy spaces}
\par  The partially ordered set $(\E_1(D),\succeq)$ is directed. Indeed, if $u,v\in\E^k_1(D)$ then $w=E(u,v)=(E(u_1,v_1),\dots,E(u_k,v_k))\in\E^k(D)$ and $w\succeq u,v$. Let $M_j$ be the total Monge--Amp\'ere mass of $w_j$. By Corollary \ref{C:mae} $M_j\le2$. Hence \[\tilde w=(M_1^{-1/n}w_1,\dots,M_k^{-1/n}w_k)\in\E^k_1(D)\] and $\tilde w\succeq u$ and $\tilde w\succeq v$.
\par By Theorem \ref{T:it} if $u\succeq v$ then $H^p_u\sbs H^p_v$ and the imbedding operator $i_{uv}$ is continuous. Thus the set of spaces $H^p_u(D)$, $u\in\E^k_1$, form a projective system (see \cite[II.6]{Sch}). Let $X^p$ be the projective limit of $(H^p_u(D),u\in\E^k_1(D))$, i. e., a subspace of all $x\in\prod_{u\in\E^k_1}H^p_u(D)$ such that $x_v=i_{uv}x_u$. Thus the mappings $i_u:\,X^p\to H^p_u(D)$ are defined. The projective  topology on $X^p$ is the weakest topology that makes all mappings $i_u$ continuous.
\par If we fix a point $w_0\in D$ and let $g(z, w_0)$ be the pluricomplex Green function with pole at $w_0$, then $u\succeq {\bf g}=(g,\dots,g)$ for all $u\in\E^k(D)$ and $H^p_u(D)\sbs H^p_{{\bf g}}(D)$. Thus we can identify $X^p$ with all $f\in H^p_{{\bf g}}(D)$ such that $f\in H^p_u(D)$ for all $u\in H^p_u(D)$.
\par Let $f$ be a holomorphic function on $D$ and $z_0\in\bd D$. The function $f$ has the admissible limit at $z_0$ if for every approach region $\A_D^\al(z_0)$ the limit
\[f^*(z_0)=\lim_{z\to z_0,z\in\A_D^\al(z_0)}f(z)\] exists.
\bT\label{T:mt} Let $f$ be a holomorphic function on a strongly pseudoconvex domain $D$ with the $C^2$ boundary. Suppose that $f$ has admissible limits at points $\{z_j\}\in\bd D$ and $\lim_{j\to\infty}f^*(z_j)=\infty$. Then for any $p>1$ there is $u\in\E_1(D)$ such that $f\not\in H^p_u(D)$.
\eT
\begin{proof} The function $\log c(z,w)$ is plurisubharmonic, negative and has a simple pole at $w$. Hence $\log c(z,w)\le g(z,w)$, where $g(z,w)$ is the pluricomplex Green function with pole at $w$. We define Green balls $G_D(w,r)=\{z\in D:\,g(z,w)\le\log r\}$. Clearly $G_D(w,r)\sbs C_D(w,r)$.
\par Let us take any positive converging series $\sum a_j$ and fix a sequence $z_j\in\bd D$ such that $f$ has admissible limits at $z_j$ and
\[\sum_{j=1}a_j^n|f^*(z_j)|^p=\infty.\]
Let $\A_j=\A_D^{\al_j}(z_j)$, where $\al_j$ are chosen so that we can find a point $w_j$ as close to $z_j$ as we want such that $G_j=G_D(w_j,e^{-1})\sbs\A_j$.
\par We will choose inductively points $w_j$. Let $w_0$ be any point. If $w_0,\dots,w_{k-1}$ have been chosen we select $w_k$ to satisfy the following conditions:
\be\item $G_k\sbs\A_k$ and $|f|>|f^*(z_k)|/2$ on $G_k$;
\item $a_jg(z,w_j)>-2^{-j-1}a_k$ on $G_k$, $0\le j\le k-1$;
\item $g(z,w_k)>-2^{-k-1}a_j$ on $G_j$, $0\le j\le k-1$.
\ee
This is possible because by Lemma \ref{L:balls} we can take $w_j$ as close to $z_j$ so that $G_k\sbs\A_k$ as we want and by \cite{C} $g(z,w)\to 0$ uniformly on compacta in $\ovr D\sm\{z_j\}$ when $w\to z_j$ and we know that $g(z,w)$ is equal to 0 on $\bd D$ when $w$ is fixed .
\par Let $u_j=a_j\max\{g(z,w_j),-2\}$. Note that if $F$ is an open set in $D$ containing $G_D(w_j,e^{-2})$ then
\[\il {F}(dd^cu_j)^n=a_j^n.\]
Let $u=E(\{u_j\})$. Since the series $v=\sum_{j=0}^\infty u_j$ converges uniformly on $\ovr D$ we see that $v\in\E$, so $u\ge v$  is a continuous plurisubharmonic function on $D$ equal to 0 on $\bd D$. Since
\[\sum_{j=0}^\infty MA(u_j)=\sum_{j=0}^\infty a_j^n<\infty\]
by Theorem \ref{T:mae} $M=MA(u)<\infty$.
\par Let us evaluate the Monge--Amp\'ere mass of $u$ on $G_k$. From the inequalities $u_k\ge u\ge v$ on $D$ and the conditions on the choices of $w_j$ on $\bd G_k$ we get
\[-a_k\ge u\ge-\sum_{j=0}^{k-1}2^{-j-1}a_k-a_k-\sum_{j=k+1}^\infty 2^{-j-1}a_k\ge-\frac32a_k.\]
Hence $u+3a_k/2\ge 0$ on $\bd G_k$ and the set $F_k=\{6(u+\frac32a_k)< u_k\}$ compactly belongs to $G_k$. Moreover if $z\in \bd G_D(w_k,e^{-2})$ then
\[6(u(z)+\frac32a_k)\le6(u_k(z)+\frac32a_k)=-3a_k<-2a_k=u_k(z).\]
Thus the set $F_k$ contains the ball $G_D(w_k,e^{-2})$. By the Comparison principle
\[6^n\il {G_k}(dd^cu)^n=\il {G_k}(dd^c6(u(z)+\frac32a_k))^n\ge6^n\il {F_k}(dd^cu_k)^n=6^na_k^n.\]
\par Hence
\[\|f\|^p_{u,p}\ge\il D|f|^p(dd^cu)^n\ge\sum_{k=0}^\infty\il {G_k}|f|^p(dd^cu)^n\ge2^{-p}\sum_{k=0}^\infty|f^*(z_k)|^pa_k^n=\infty.\]
Hence $f\not\in H^p_u(D)$.
\end{proof}
\par Let us introduce a new topology on the space $(H^\infty(D))^k$. Consider imbeddings $j_u:\,(H^\infty(D))^k\to H^p_u(D)$, $u\in\E_1^k$, and for any $R>0$ the sets $j_u^{-1}(B_{u,p}(R))$. These sets with an empty set form a basis because they, evidently, cover $(H^\infty(D))^k$ and for any $x,y\in (H^\infty(D)^k$, any $u,v\in\E^k_u(D)$ and any $R_1,R_2>0$ the intersection $A$ of the sets $x+i_u^{-1}(B_{u,p}(R_1))$ and $y+i_v^{-1}(B_{v,p}(R_2))$ contains an element of the basis. Indeed, if $A$ is empty then there is nothing to prove. If $z\in A$ then $\|z-x\|_{u,p}<R_1$ and $\|z-y\|_{v,p}<R_2$. Let $w=E(u,v)$ and $\tilde w=(\al_1w_1,\dots,\al_kw_k)$, where the coefficients $\al_j\ge 1/2$ have been chosen so that $\tilde w\in\E^k_1(D)$. Since $\|f\|_{w,p}\ge\max\{\|f\|_{u,p},\|f\|_{v,p}\}$ we see that $\|f\|_{\tilde w,p}\ge2^{-n/p}\max\{\|f\|_{u,p},\|f\|_{v,p}\}$. Hence $B_{\tilde w,p}(2^{-n/p}R)\sbs B_{u,p}(R)\cap B_{v,p}(R)$ and we see that there is $c>0$ such that the set $z+B_{\tilde w,p}(c)\sbs A$.
\par We denote by $Y^p$ the space $(H^\infty(D))^k$ endowed with the topology defined by the basis of sets $j_u^{-1}(B_{u,p}(R))$ for all $u\in\E_1^k$ and all $R>0$.
\bT\label{T:int} Let $D$ be a strongly pseudoconvex domain with the $C^2$ boundary and let $p\ge1$. Then $\cap_{u\in\E^k_1(D)}H^p_u(D)=(H^\infty(D))^k$ and the projective limit $X^p$ of $(H^p_u(D),u\in\E^k_1(D))$ is isomorphic to $Y^p$.
\eT
\begin{proof} It suffices to prove this theorem for $k=1$. Since all mappings $i_{uv}$ are imbeddings if $x\in X^p$ and $x=(f_u,u\in\E_1)$ then $f_u=f_v=f$ and this $f$ belongs to all spaces $H^p_u$ or $f\in \cap_{u\in\E^k_1(D)}H^p_u(D)$. Let us show that the latter space is $(H^\infty(D))^k$. Suppose that $f$ be unbounded. Since $f\in H^p_{{\bf g}}$ by \cite[Theorem 10]{S} $f$ has admissible limits a.e. on the boundary.  If the function $f^*$ is bounded then the real and imaginary parts of $f$, which are harmonic functions, have bounded admissible limits equal to $f^*$ a.e. (see \cite{AS, Sm}) and this implies that $f$ is bounded. Hence $f^*$ is unbounded. By Theorem \ref{T:mt} there is $u\in\E_1(D)$ such that $f\not\in H^p_u(D)$. Thus $f\in H^\infty(D)$ and we got a mapping $\Phi:\,X^p\to Y^p$. Clearly, this mapping is an algebraic isomorphism.
\par  By its definition the projective topology on $X^p$ must contain all sets $A(x,u,R)=x+i_u^{-1}(B_{u,p}(R))$, where $u\in\E^k_1(D)$, $x\in (H^\infty(D)^k$ and $R>0$. It is easy to see that $F(A(x,u,R))=\Phi(x)+j_u^{-1}(B_{u,p}(R))$. We conclude that these sets form a basis of the projective topology on $X^p$ and, therefore, $\Phi$ is a topological isomorphism.
\end{proof}
\par The duals of $H^p_u(D)$ form an inductive system and their inductive limit can be considered. We will not go here into this. Instead, we will show that the intersection of any countable family of spaces $H^p_u(D)$ contains an unbounded function.
\bT Let $D$ be a strongly pseudoconvex domain with the $C^2$ boundary and let $p\ge1$. Let $\{u_j\}\sbs\E_1(D)$. Then the space $X=\cap_{j=1}^\infty H^p_{u_j}(D)$ contains an unbounded function.
\eT
\begin{proof} Let us pick up positive coefficients $\al_j$ such that the function \[u=\sum_{j=1}^\infty \al_ju_j\in\E_1(D).\] Clearly $H^p_u(D)\sbs X$ and we only need to prove that for any $u\in\E_1(D)$ the space $H^p_u(D)$ contains an unbounded function. If not then the continuous imbedding $H^\infty(D)\to H^p_u(D)$ is onto. By a theorem of Banach the inverse mapping is also continuous. Let us find a point $z_0\in\bd D$ such that $\mu_u(\{z_0\})=0$ and take a peak function $q$ at $z_0$. The norm of the functions $q^m$, $m\in\aN$, in $H^\infty(D)$ is 1.  \cite[Theorem 3.1]{D} states that for a plurisubharmonic function $\phi$ on $D$ continuous up to the boundary
\[\mu_u(\phi)=\il {B_u(r)}\phi(dd^cu)^n+\il
{B_u(r)}(r-u)dd^c\phi\wedge(dd^cu)^{n-1}.\] Thus the norms of the functions $q^m$ in $H^p_u(D)$ are equal to $\mu_u^{1/p}(|q|^{pm})$ and, consequently, converge to 0. We came to a contradiction.
\end{proof}

\end{document}